
\input amstex

\documentstyle{amsppt}
\NoBlackBoxes

\define\scrO{\Cal O}
\define\Pee{{\Bbb P}}
\define\Zee{{\Bbb Z}}
\define\Cee{{\Bbb C}}
\define\Ar{{\Bbb R}}

\define\dirac{\rlap{/}\partial}
\define\dbar{\bar\partial}
\define\Spin{\operatorname{Spin}}
\define\Id{\operatorname{Id}}

\topmatter
\title Donaldson and Seiberg-Witten Invariants of Algebraic Surfaces\endtitle
\author Robert Friedman\endauthor

\address Department of Mathematics, Columbia University, New York NY 10027
\endaddress

\email rf\@math.columbia.edu\endemail


\thanks The author was supported in part by NSF
Grant \# DMS-92-03940.\endthanks

\endtopmatter

\document

\head 1. Introduction \endhead

Donaldson theory and more recently Seiberg-Witten theory have led to dramatic
breakthroughs in the study of smooth $4$-manifolds, and in particular of
algebraic surfaces and their generalizations, symplectic $4$-manifolds. In 
this paper, we shall survey some of the main results. Many expositions of
Seiberg-Witten theory have appeared, for example the survey article of
Donaldson \cite{10}  and the book of Morgan 
\cite{28}, and the only claim to originality made here is
in the focus on algebraic surfaces. A general reference for Donaldson theory
is the book of Donaldson and Kronheimer \cite{11}. A discussion of the
invariants from the point of view of physics has been given by Witten
\cite{41}. A more detailed  exposition of the results described here for 
K\"ahler manifolds is given in the papers of the author and Morgan 
\cite{16} and Brussee
\cite{5}, as well as in  Okonek-Teleman \cite{33}, \cite{34}.

  From the viewpoint of algebraic geometry, the major numerical invariants of   
an algebraic surface $X$ are the plurigenera $P_n(X) = H^0(X; nK_X)$,     
$n\geq 1$. Indeed, the second plurigenus $P_2(X)$ already appears in
Castelnuovo's criterion for rationality, and invariants $P_4(X)$ and    
$P_6(X)$ appear in Enriques' criterion for when a surface is rational or  
ruled. Unlike the more natural invariants
$p_g(X)$ and $q(X)$, however, which are linked to the oriented homotopy   
theory of $X$ via Hodge theory, there is no natural topological   
interpretation (as yet) for the higher invariants $P_n(X)$ for $n\geq 2$.      
A deep geometric fact about surfaces is the existence of strong minimal   
models for algebraic surfaces of Kodaira dimension at least zero (in other
words, those surfaces which are not rational or ruled, or equivalently for 
which $P_n(X) \neq 0$ for some $n$).  Assuming for simplicity that $X$ is the
blowup of a minimal surface $X_0$ at distinct points, the curve fibers of the
blowup morphism are then distinguished  holomorphically embedded copies of 
$\Pee ^1$ in
$X$ with self-intersection $-1$, and their cohomology classes are special
cohomology classes as far as the algebraic geometry of $X$ is concerned. A
further marked class is that of
$K_X$. However, from the point of view of diffeomorphism, it is better to   
work with the pullback of $K_{X_0}$, since if $X$ is not minimal then 
$c_1(K_X)$ is almost never invariant up to sign under diffeomorphisms of $X$.

The main questions in 1935 concerning the topology of algebraic surfaces
described in, say, Zariski's book \cite{42}  are mainly concerned with what
would now be viewed as the consequences of Hodge theory, and as such give
invariants which depend only on the oriented homotopy type of $X$. Questions
more directly concerned with the smooth topology of $X$ date back to Severi's
problem \cite{38} of giving a topological criterion for rationality (1949).
Actually, given the examples of the Enriques and Godeaux surfaces, Severi  
asked if a surface with
$p_g = 0$ and
$H_1(X; \Zee) = 0$ was necessarily rational. Later Dolgachev \cite{6} 
constructed examples of simply connected elliptic surfaces which were not
rational (1967), and Barlow
\cite{1}  gave the only known example of a simply connected surface of
general type with
$p_g=0$ (1980). By Freedman's work (1981), it is known that these surfaces   
are all homeomorphic to rational surfaces. In a related vein, Kodaira
\cite{21} classified all of the complex surfaces homotopy equivalent to $K3$
surfaces (1967) and asked if they were homeomorphic to $K3$ surfaces; again 
this is a consequence of Freedman's theorem. Thanks to the examples of 
Dolgachev and Kodaira, elliptic surfaces were intensively investigated    
during the 1970's by Moishezon and Mandelbaum (see for example \cite{27}),    
as well as by Harer-Kas-Kirby \cite{19}, using Kirby's
handlebody calculus. Using these techniques one can show for example that
all homotopy equivalent simply connected elliptic surfaces become  
diffeomorphic once we connect sum one copy of $\Cee P^2$. However, these
techniques proved unable to give much information about the actual
diffeomorphism classification of elliptic surfaces, and matters remained at   
an impasse until Donaldson's famous counterexample \cite{8} to the
$h$-cobordism theorem in dimension $4$, which showed via gauge theory that     
a certain Dolgachev surface was not diffeomorphic to a rational surface. At 
this point, the power of gauge theory as a tool in attacking the smooth
classification of complex surfaces became apparent, and it was natural to
conjecture that the smooth topology of an algebraic surface would in many   
ways closely reflect its algebro-geometric structure. A related question is  
the Thom conjecture and its generalizations. The following  is a list of 
theorems along these lines which can be proved  by Donaldson theory or by
Seiberg-Witten theory (which has also simplified the previous proofs in case 
the result was already known by Donaldson theory).

\proclaim{Theorem} For a fixed diffeomorphism type of a $4$-manifold $M$,  
there are only finitely many complex structures on $M$ up to deformation of
complex structure, at least so that the resulting complex surface is   
K\"ahler.
\endproclaim
A proof is given in the book \cite{15}.
(In the non-K\"ahler case, the difficulty lies in the classification of 
surfaces of type VII in Kodaira's notation, but presumably the result is still
true.) The main part of the argument involves studying elliptic surfaces, and
these surfaces have now been classified up to diffeomorphism \cite{13},
\cite{15},
\cite{31}, \cite{2}, \cite{3}, \cite{12}, \cite{30}:

\proclaim{Theorem} Two elliptic surfaces with positive holomorphic Euler
characteristic are diffeomorphic if and only if they are deformation 
equivalent.
\endproclaim
As a consequence, there is an answer to the question raised by Kodaira: 

\proclaim{Theorem} A complex
surface diffeomorphic to a
$K3$ surface is a $K3$ surface.
\endproclaim

One can also classify elliptic surfaces with holomorphic Euler
characteristic equal to zero, and show that diffeomorphism implies   
deformation equivalence up to a two-to-one ambiguity caused by orientation
questions
\cite{15}. Thus, for elliptic surfaces $X$ and $X'$, if $X$ is diffeomorphic  
to
$X'$ then $X$ is deformation equivalent either to $X'$ or to its complex
conjugate.

\proclaim{Theorem (the van de Ven conjecture)} If $X$ and $X'$ are 
diffeomorphic complex surfaces, then $X$ and $X'$ have the same Kodaira
dimension.
\endproclaim 

 A proof is given in  \cite{18} as well as
\cite{36} and \cite{35}. A main step in the proof is the answer to the  
question raised by Severi: 

\proclaim{Theorem} A complex surface diffeomorphic to a rational surface is
rational. 
\endproclaim

Using Seiberg-Witten theory, one can prove a stronger statement: 

\proclaim{Theorem} If
$X$ and $X'$ are diffeomorphic complex surfaces, then $P_n(X) = P_n(X')$ for
all $n\geq 1$. 
\endproclaim

For an elliptic surface, $P_n(X)$ can be calculated from the homotopy type
once we know the multiple fibers, so the theorem for elliptic surfaces
follows from the diffeomorphism classification of elliptic surfaces. For a
surface $X$ of general type, the plurigenera are determined by $K_{X_0}^2$,
where $X_0$ is the minimal model of $X$. Thus the  proof of this theorem 
follows once we know that the $C^\infty$ topology of
$X$ determines the number of times $X$ is blown up from its minimal
model. 

\proclaim{Theorem} Suppose that $X$ and $X'$ are complex surfaces with   
$\kappa (X) \geq 0$, and that $f\: X \to X'$ is an orientation-preserving
diffeomorphism. \rom(Thus necessarily $\kappa (X') \geq 0$ as well.\rom) Let
$X_0$ be the minimal model of
$X$ and $X_0'$ the minimal model of $X'$. Likewise let $[E_i], 1\leq i      
\leq r$, be the classes of the exceptional curves on $X$ and let $[E_j'],  
1\leq j\leq s$, be the classes of the exceptional curves on $X$. Then $r=s$  
and for all $j$, $f^*[E_j'] = \pm [E_i]$ for some $i$. Moreover, after
identifying
$[K_{X_0}]$ with a cohomology class in $H^2(X; \Zee)$ and likewise for
$[K_{X_0'}]$, we have $f^*[K_{X_0'}] = \pm [K_{X_0}]$. 
\endproclaim

This theorem, whose proof
seemed inaccessible to Donaldson theory, follows quite easily from the form   
of the Seiberg-Witten invariants for a K\"ahler surface as worked out by 
Witten, Taubes, and Kronheimer-Mrowka, and seems to have been first   
explicitly noted for K\"ahler surfaces with $b_2^+\geq 3$ by Kronheimer. The
case $b_2^+ = 1$, or equivalently $p_g=0$, needs some further analysis of
chamber structures which is rather easy to handle in the Seiberg-Witten case 
and is worked out in \cite{16} as well as \cite{5}. 

\proclaim{Theorem} Suppose that $X$ is a complex surface and that $X$ is
diffeomorphic to a connected sum $M_1\#M_2$. Then one of
$M_1$, $M_2$ is negative definite \rom(Donaldson \cite{9}\rom). If       
$\kappa (X)
\geq 0$ and, say, $M_2$ is
negative definite, then $H_2(M_2)$ is contained in the subgroup of
$H_2(X)$ spanned by the classes of the exceptional curves on $X$. In
particular, if  $\kappa (X) \geq 0$, then every smoothly embedded $2$-sphere  
in
$X$ with self-intersection $-1$ is homologous up to sign to the class of an
exceptional curve.
\endproclaim

Again this
follows in a straightforward way in case $b_2^+\geq 3$, and for
$b_2^+ = 1$ by analyzing the chamber structure.

Finally there is the Thom conjecture and its generalizations: 

\proclaim{Theorem} Let $X$ be a
K\"ahler surface and let $\alpha \in H_2(X;\Zee)$ satisfy $\alpha ^2 \geq 0$.
If
$\alpha$ is represented by a closed oriented smooth $2$-manifold of genus $g$,
then $2g-2 \geq \alpha ^2 + \alpha \cdot K_X$. 
\endproclaim
Of course, by adjunction, if
$\alpha = [C]$ where $C$ is a holomorphically embedded algebraic curve, then
$2g(C) - 2 = \alpha ^2 + \alpha \cdot K_X$, and so the genus of $C$ is minimal
among all closed oriented smooth $2$-manifolds whose fundamental class is
$\alpha$. The Thom conjecture is the special case $X = \Pee ^2$. Following a
program outlined by Kronheimer \cite{22}, Kronheimer and Mrowka attacked this
conjecture in a series of deep papers \cite{23} and proved many special cases
(for example in case $X$ is a $K3$ surface). In particular they were led to  
the discovery of the structure of Donaldson polynomials for manifolds of  
simple type (to be explained below). The general conjecture was proved via
Seiberg-Witten theory by Kronheimer and Mrowka \cite{26} and independently by
Morgan-Szab\'o-Taubes
\cite{32}. 

\head 2. A brief review of Donaldson theory \endhead

As we shall see, one fundamental reason that Donaldson theory and
Seiberg-Witten theory exist to give such powerful information about dimension
$4$ is that $4$ is the unique solution to the equation $2= n-2$.

Let $M$ be a closed oriented Riemannian $4$-manifold and let $P\to M$ be a
principal $SU(2)$-bundle over $M$ (the theory also studies the case of an 
$SO(3)$-bundle, but we shall not discuss this). Such bundles are classified   
by the integer $c = c_2(P) \in H^4(M; \Zee) \cong \Zee$. If $A$ is a  
connection on
$P$, then its curvature $F_A$ is a $2$-form with values in $\operatorname{ad}
P$. It is natural to consider the energy
$$\int _M|F_A|^2,  \qquad \text{the {\sl Yang-Mills functional}},$$
and to try to minimize this energy. In this sense, Yang-Mills theory is a
nonabelian analogue of Hodge theory, which attempts to minimize the norm of a
form representing a fixed cohomology class. The critical points of the
Yang-Mills functional are given by the Euler-Lagrange equations $D_AF_A = 0,
D_A*F_A = 0$. Here $D_A$ is the differential operator on sections of $\Omega
^2(M; \operatorname{ad} P)$ associated to $A$,  the equation  $D_AF_A = 0$
is the Bianchi identity which is always satisfied, and  $*$ is the Hodge
$*$-operator from $\Omega ^2(M; \operatorname{ad} P)$ to $\Omega ^{n-2}(M;
\operatorname{ad} P)$. In case
$\dim M = 4$,
$*F_A$ is also a section of $\Omega ^2(M; \operatorname{ad} P)$. Thus if  
$*F_A$ is a scalar multiple of $F_A$, then the Bianchi identity also implies 
that
$D_A*F_A = 0$. Since $** = \operatorname{Id}$, the only possibilities are  
$*F_A = F_A$, in which case we say that $A$ is a {\sl self-dual\/} connection, 
or
$*F_A =- F_A$,  in which case we say that $A$ is an {\sl anti-self-dual\/} 
(ASD) connection. These connections, if they exist, are absolute minima for  
the Yang-Mills functional. 

For a $4$-manifold $M$,  we can break up the space of sections $\Omega ^2(M;
\operatorname{ad} P)$ into its $+1$ and $-1$ eigenspaces for the Hodge
$*$-operator: 
$$\Omega ^2(M; \operatorname{ad} P) = \Omega ^2_+(M;
\operatorname{ad} P) \oplus \Omega ^2_-(M; \operatorname{ad} P).$$
Write $F_A^+$ for the component of $F_A$ lying in $\Omega ^2_+(M;
\operatorname{ad} P)$, so that $A$ is ASD if and only if $F_A^+ = 0$ if and 
only if $F_A \in \Omega ^2_-(M; \operatorname{ad} P)$. Given the bundle $P$  
and the metric $g$, the set of all anti-self-dual connections, modulo the  
group of symmetries of $P$, forms a finite-dimensional oriented smooth  
manifold
$\Cal M(P,g)$, at least for a generic metric $g$. While
$\Cal M(P,g)$ is not in general compact, it can be compactified to a  
stratified space $X(P, g)$ which carries a fundamental class $[X(P, g)]$.    
Now morally speaking there is a tautological $SU(2)$-bundle $\Cal P$ over
$M\times
\Cal M(P,g)$. Taking the slant product with $c_2(\Cal P)$ induces a 
homomorphism
$\mu \: H_*(M; \Zee) \to H^{4-*}(\Cal M(P,g); \Zee)$. Again morally speaking
the classes $\mu (\alpha)$ extend to classes in $H^{4-*}(X(P, g); \Zee)$,    
and appropriate combinations can then be evaluated on the fundamental class
$[X(P, g)]$ to produce polynomial invariants for $M$, which will be  
independent of $g$ as long as $b_2^+(M) >1$.  

Why should Donaldson theory give more information about the smooth topology 
than the standard homotopy or homeomorphism invariants? This is a deep  
mystery, but one could attempt to give a partial answer, or perhaps an
ideological underpinning to Donaldson theory, by making the following points
about the structure of the ASD equations and the corresponding moduli spaces:
\roster
\item The groups $SU(2)$ and $SO(3)$ are nonabelian, and so the associated
nonabelian Hodge theory might be able to reveal more structure than the usual
abelian Hodge theory, which only gives homotopy-theoretic information about 
$M$.
\item The ASD equations are conformally invariant.
\item As a consequence of the conformal invariance, which allows for bubbling
off, the ASD moduli spaces are almost never compact, and the compactification
involves the underlying manifold $M$ in an interesting and nontrivial way.
(Sometimes, however, even compact moduli spaces can contain the essential
information in Donaldson theory.)
\item The theory works best for simply connected $4$-manifolds. One can define
the polynomial invariants for arbitrary $4$-manifolds, but the definition is
much more difficult \cite{29}, and there are ideological grounds for
concentrating on the class of simply connected $4$-manifolds. (Arbitrary
$4$-manifolds are technically unknowable, because every finitely presented
group is the fundamental group of a smooth $4$-manifold, and on the other hand,
modulo the three-dimensional Poincar\'e conjecture, dimension $4$ is the first
dimension where we encounter really interesting simply connected manifolds.)
\endroster

Now suppose that $M$ is a K\"ahler surface $X$, with K\"ahler form $\omega   
\in
\Omega ^{1,1}(X)$. Then it is an easy calculation using the Hodge identities
that
$$\align
\Omega ^2_+(X; \operatorname{ad} P)\otimes \Cee &= \Omega ^{0,2}(X;
\operatorname{ad} P) \oplus \Omega ^{2,0}(X; \operatorname{ad} P)\oplus \Omega
^0(X; \Cee) \cdot \omega;\\
\Omega ^2_-(X; \operatorname{ad} P)\otimes \Cee &= (\omega)^\perp \cap \Omega
^{1,1}(X; \operatorname{ad} P),
\endalign$$
where this last equation means that the elements in $\Omega ^2_-(X;
\operatorname{ad} P)\otimes \Cee$ are $(1,1)$-forms pointwise orthogonal to
$\omega$. Let
$F_A^{0,2}$ be the projection of the
$2$-form $F_A$ onto the subspace $\Omega ^{0,2}(X; \operatorname{ad} P)$, and
similarly for $F_A^{1,1}$. Thus the equation $F_A^+ = 0$ is equivalent to the
two equations
$$\align
F_A^{0,2} &= 0;\\
F_A^{1,1} &\text{\, \,is pointwise orthogonal to $\omega$}.
\endalign$$
We can also write this last equation as $\Lambda F_A^{1,1} = 0$, where the
operator $\Lambda$ is contraction against $\omega$. The first equation says
that the projection $A^{0,1}$ of $A$ onto the $(0,1)$-forms defines the
$\dbar$-operator for a holomorphic structure, with holomorphically trivial
determinant, on the complex $2$-plane bundle associated to $P$. This is an  
easy integrability result, which follows for example from the 
Newlander-Nirenberg theorem. The second equation, according to the
Kobayashi-Hitchin conjecture, proved by Donaldson \cite{7} in the case of
K\"ahler surfaces and by Uhlenbeck-Yau \cite{40} for general K\"ahler 
manifolds, is equivalent to saying that the holomorphic vector bundle  
structure on $V$ is
$\omega$-stable. In the case of a rank two holomorphic vector bundle $V$ on a
surface $X$ with $c_1(V) = 0$, this condition amounts to the following: if $L$
is a line bundle and there is a nonzero holomorphic map
$L
\to V$ (not necessarily an inclusion on each fiber), then $c_1(L) \cdot  
\omega < 0$. In case $\omega$ is the K\"ahler form of a Hodge metric
corresponding to an ample divisor $H$, then $\omega$-stability exactly
corresponds to stability with respect to $H$, a notion introduced by Mumford 
and later Takemoto in order to construct moduli spaces on curves and
higher-dimensional projective varieties. While not completely successful in
constructing compact moduli spaces (one needs instead to consider all   
Gieseker semistable torsion free sheaves), stability is a fundamental
non-degeneracy condition for vector bundles on a projective variety. As such, 
it has  numerical consequences for $V$, in the form of Bogomolov's inequality. 
In the case of $V$ of rank two and trivial determinant, this inequality simply
reads $c_2(V) \geq 0$, and is an easy consequence of the existence of an ASD
metric on $V$. Moreover, again by using the existence of an ASD metric, the 
case
$c_2(V) = 0$ corresponds to the case where $V$ is flat, in other words
associated to an irreducible representation of $\pi _1(X)$ in $SU(2)$. But   
not every bundle which satisfies the numerical condition $c_2\geq 0$ is
automatically stable, so that there are also more subtle, not strictly 
numerical consequences of stability which are necessary in order to    
construct separated moduli spaces.

In any case, the connection between ASD connections and stable holomorphic
vector bundles leads to some calculation of polynomial invariants by
algebro-geometric methods, and suggests that there is a link between the
algebraic geometry of a K\"ahler surface and its $4$-manifold topology. Of
course, a little experience shows that the moduli spaces of vector bundles on
algebraic surfaces are very subtle invariants to calculate. Already for $\Pee
^2$, they become quite complicated when $c_2$ grows, and as the complexity of
the surface increases, the complexity of the corresponding moduli spaces also
seems to increase.

To return to the general story of Donaldson theory, one can ask if the
Donaldson polynomial invariants can detect special cohomology classes on
$M$ (where $b_2^+(M) \geq 3$). One way to do this is the following: under 
certain circumstances, the Donaldson polynomials lie in the subring   
$\Cee[q_M, k]$ of the symmetric algebra, where
$q_M$ is the symmetric polynomial corresponding to the intersection form on  
$M$ and $k$ is a class in $H^2(M; \Zee)$. If there exists a polynomial of   
this type not lying in $\Cee[q_M]$, then it is an easy algebraic argument   
that the class $k$ is preserved up to $\pm 1$ by
orientation-preserving self-diffeomorphisms of
$M$. Another circumstance where Donaldson theory can find special classes
is the following: suppose that there is a smoothly embedded $2$-sphere in
$M$ whose associated cohomology class $\alpha$ satisfies $\alpha ^2 = -1$;
for example this is the case for the classes of exceptional curves. Then
the Donaldson polynomials have a rather restricted form. If there
exists a nonzero Donaldson invariant for $M$, then the possible such classes
$\alpha$ represented by smoothly embedded $2$-spheres of self-intersection
$-1$ are either equal up to $\pm 1$ or orthogonal, and in particular there
are only finitely many such classes. For a description of these and
related results, see \cite{15}. After a period involving a considerable amount 
of difficult calculation of Donaldson polynomials, by algebro-geometric and 
other means, Kronheimer and Mrowka proved a deep structure theorem for a large
class of
$4$-manifolds, those of {\sl simple type}. While we shall not give a precise
definition, a
$4$-manifold $M$ with $b_2^+(M) \geq 3$ is of simple type if evaluating the
polynomial invariant on the $4$-dimensional  class $\mu (\text{pt})$ 
essentially just gives back a polynomial invariant of a smaller-dimensional
moduli space. Large classes  of $4$-manifolds have simple type, and there is  
no known example of a
$4$-manifold with $b_2^+\geq 3$ which does not have simple type. Kronheimer
and Mrowka \cite{24, 25} showed that, if $M$ is a simply connected  
$4$-manifold of simple type, then  the Donaldson polynomials can be expressed
via certain recurrence relations as polynomials involving certain {\sl basic
classes\/}
$\kappa _i \in H^2(M; \Zee)$ and rational numbers $a_i$. Here the $\kappa _i$
are characteristic elements of $H^2(M; \Zee)$, i\.e\. for all $\alpha \in
H^2(M; \Zee)$, $\alpha ^2 \equiv \alpha \cdot \kappa \bmod 2$.  

\head 3. Definition of the Seiberg-Witten invariants \endhead

One basic motivation which led to the definition of the Seiberg-Witten
invariants was the attempt to find an {\it a priori\/} description of the
Kronheimer-Mrowka basic classes. To define the Seiberg-Witten invariants, we 
need to recall the definition of a
$\Spin ^c$ structure on $M$. Let $M$ be an oriented Riemannian $4$-manifold
(actually we can define a $\Spin ^c$ structure in every dimension). Then the
tangent bundle of $M$ corresponds to a principal $SO(4)$-bundle. Now $SO(4)$ 
has a unique double cover $\Spin (4)$, which by one of the coincidences of Lie
group theory in low dimensions is isomorphic to $SU(2) \times SU(2)$, and one
can ask if the tangent bundle of
$M$ lifts to a $\Spin (4)$ bundle. The answer is that such a lift exists if  
and only if the second Stiefel-Whitney class $w_2(M)$ is zero, and in case $M$ 
is simply connected (or more generally if $H^2(M; \Zee)$ has no $2$-torsion) 
this condition is equivalent to assuming that
$H^2(M; \Zee)$ has an even intersection form. In this case a lift of $TM$ to
$\Spin (4)$ is called a {\sl $\Spin$ structure\/} on $M$. Of course, most
$4$-manifolds do not have a $\Spin$ structure. However, if one does exist,  
then from the isomorphism $\Spin (4) \cong SU(2) \times SU(2)$, there are two
associated rank two complex vector bundles $\Bbb S^+, \Bbb S^-$ with trivial
determinant. A fundamental fact is that the Levi-Civita connection on $TM$
induces a differential operator, the {\sl Dirac operator\/} $\dirac\: \Bbb S^+
\to
\Bbb S^-$. Here
$\dirac$ is a first order elliptic formally self-adjoint operator.

If there is no $\Spin$ structure, there is still a related construction.  
Define
$$\Spin ^c(4) = \Spin (4) \times U(1)/\{\pm 1\},$$
where the group $\{\pm 1\}$ acts diagonally on both factors. While these  
groups can be defined in general, in dimension $4$ we also have the  
isomorphism
$$\Spin ^c(4) \cong \{\, (A,B) \in U(2) \times U(2): \det A = \det B\,\}.$$
By definition a $\Spin ^c$ structure on $M$ is a lift of the $SO(4)$-bundle 
$TM$ to a bundle with structure group $\Spin ^c(4)$. 
There is an exact sequence
$$\{1\} \to U(1) \to \Spin ^c(4) \to SO(4) \to \{1\},$$
and from this exact sequence it is straightforward to see that $M$ has a  
$\Spin ^c$ structure if and only if the image of $w_2(M)$ under the Bockstein
homomorphism $H^2(M; \Zee/2\Zee) \to H^3(M; \Zee)$ is zero, or equivalently   
if $w_2(M)$ is the mod $2$ reduction of an integral class. In fact, this  
condition is satisfied by every $4$-manifold $M$. In case $M$ is an almost
complex surface
$X$, a natural lift of $w_2(X)$ is $c_1(X)$. The set of all $\Spin ^c$ 
structures on $M$ is a principal homogeneous space over $H^2(M; \Zee)$.

Let $\xi$ be a $\Spin ^c$ structure on $M$. Since $\Spin ^c(4) \subset U(2)
\times U(2)$, there are two associated complex $2$-plane bundles $\Bbb
S^+(\xi), \Bbb S^-(\xi)$, and $\det  \Bbb S^+(\xi) = \det \Bbb S^-(\xi)$. We
call the complex line bundle $L = \det  \Bbb S^+(\xi) = \det \Bbb S^-(\xi)$ 
the  {\sl determinant\/} of $\xi$ and will write it as $\det \xi$. Note that  
if
$\xi$ and $\xi '$ differ by $\alpha \in H^2(M; \Zee)$, then $\Bbb
S^{\pm}(\xi') = \Bbb S^{\pm}(\xi) \otimes L_\alpha$, where $L_\alpha$ is
a complex line bundle such that $c_1(L_\alpha ) = \alpha$, and $c_1(\det
\xi ') = c_1(\det \xi) + 2\alpha$. In particular, the mod two reduction of
$c_1(\det \xi)$ depends only on $M$, and it is easy to see that this reduction 
is exactly $w_2(M)$. Choosing a connection $A$ on $\det \xi$ gives a Dirac
operator $\dirac _A\: \Bbb S^+(\xi) \to \Bbb S^-(\xi)$, which is again a first
order elliptic formally self-adjoint operator. The index of $\dirac _A$ is
given by the Atiyah-Singer index theorem:
$$\operatorname{index} \dirac _A = \frac18(c_1(\det \xi)^2 - \sigma (M)),$$
where $\sigma (M) = b_2^+(M) - b_2^-(M)$ is the signature of $M$. 

Now suppose that $M=X$ is an almost complex $4$-manifold. Then the structure
group of $M$ reduces to $U(2)$, and we seek a lift of the inclusion
$U(2)\subset  SO(4)$ to the group $\Spin ^c(4) \subset U(2) \times U(2)$. One
checks that a natural lift is given by
$$T \in U(2) \mapsto \left(\pmatrix \Id &0\\0 & \det T\endpmatrix, T\right)  
\in U(2) \times U(2).$$ 
For this lift, if $\xi _0$ denotes the corresponding $\Spin ^c$ structure, we
have $\Bbb S^+(\xi_0) = \Omega ^0(X) \oplus K_X^{-1}=
\Omega ^0(X) \oplus \Omega ^{0,2}(X)$ and $\Bbb S^-(\xi_0) = TX= \Omega
^{0,1}(X)$, where $K_X$ is the canonical bundle of the almost complex 
structure. Thus every
$\Spin ^c$ structure $\xi$ on $X$ is given by choosing a complex line bundle
$L_0$ on
$M$, corresponding to a class in $H^2(X; \Zee)$, and twisting $\xi _0$ and 
$\Bbb S^{\pm}(\xi_0)$ by $L_0$. In this case
$$\Bbb S^+(\xi) = \Omega ^0(L_0) \oplus \Omega ^{0,2}(L_0), \qquad          
\Bbb S^-(\xi) =  \Omega ^{0,1}(L_0),$$
and $\det \xi = L_0^2\otimes K_X^{-1}$, or equivalently $L_0$ is a square root
of $K_X\otimes L$, where $L = \det \xi$ is a characteristic line bundle (the
mod two reduction of $c_1(L)$ is $w_2(X)$). There are two obvious
choices (not necessarily distinct) for a $\Spin ^c$ structure on $X$: we
can choose $L_0$ to be trivial and $L = K_X^{-1}$ or $L_0 = K_X$ and
$L=K_X$. We will refer to these two choices as the {\sl trivial\/} $\Spin ^c$
structures on $X$. Finally, in case
$X$ is K\"ahler, the operator $\dirac _A$ corresponds under this isomorphism  
to
$\sqrt{2}(\dbar _A+ \dbar _A^*)$. Here $\dbar _A$ is the $\dbar$-operator on
$X$ coupled with the connection induced on $L_0$ by $A$, which can be
defined for every $C^\infty$ complex line bundle $L$ (not necessarily
holomorphic) and connection $A$ on $L$. This is essentially a calculation in
Euclidean space, based on the fact that a Hermitian metric is K\"ahler if and
only if it looks like the standard metric to second order.

We can now give the Seiberg-Witten equations for a smooth
Riemannian $4$-manifold $M$: let $\xi$ be a $\Spin ^c$ structure on $M$, with
$\det \xi = L$. Then the equations, for a $C^\infty$ section $\psi$ of $\Bbb
S^+(\xi)$ and a connection $A$ on $L$ are as follows:
$$\align
\dirac_A\psi &= 0;\\
F_A^+ &= \psi \otimes \psi ^* - \frac{|\psi|}{2}^2\Id.
\endalign$$
Here the first equation can be paraphrased by saying that $\psi$ is a   
harmonic spinor. The second equation means the following: using the metric    
to identify 
$\Bbb S^+(\xi)$ with its dual, both $\psi \otimes \psi ^*$ and
$\dsize\frac{|\psi|}{2}^2\Id$ are sections of $\Bbb S^+(\xi) \otimes \Bbb
S^+(\xi)^* = \operatorname{Hom}(\Bbb S^+(\xi), \Bbb S^+(\xi))$. The   
difference has trace zero. Now using Clifford multiplication on   $\Bbb
S^+(\xi)$, one can identify the traceless endomorphisms from $\Bbb S^+(\xi)$  
to itself with
$\Omega ^2_+(M) \otimes \Cee$, in such a way that $\dsize\psi \otimes \psi ^*  
- \frac{|\psi|}{2}^2\Id$ becomes a purely imaginary self-dual $2$-form, and we
can compare it with $F_A^+$. There is also a symmetry group, the gauge group  
of automorphisms of $\xi$ (which are understood to induce the identity on the
tangent bundle). Up to a factor of $2$, this group is just the group of
automorphisms of the line bundle $L$ and hence it is abelian. The quotient
space of the set of solutions to the Seiberg-Witten equations by the gauge
group is the {\sl Seiberg-Witten moduli space}. It is a compact space locally
modeled on a real analytic variety. A solution $(\psi, A)$ of the
Seiberg-Witten equations is {\sl reducible\/} if $\psi = 0$, in which
case $A$ is necessarily a flat connection. For simplicity we shall
assume henceforth that there are no reducible solutions to the Seiberg-Witten
equations. Provided that
$b_2^+(M) > 0$, for generic perturbations of the equations, the   
Seiberg-Witten moduli space becomes a smooth manifold $\Cal M(\xi)$ of 
dimension 
$$\frac14\left(L^2 - (2\chi (M) + 3\sigma (M)\right),$$
and in particular  $\Cal M(\xi)$ is empty whenever this integer is negative.  
An orientation on  $\Cal M(\xi)$ can be given by choosing orientations on 
$H^0(M;
\Ar), H^1(M; \Ar)$, and $H^2_+(M; \Ar)$. Thus for example once we have    
chosen such an orientation, in case  $\Cal M(\xi)$ has dimension zero it is
a signed collection of points, and we can add up these signs to obtain an
integer. Changing the orientation simply changes the sign of this integer. In
case $\Cal M(\xi)$ has dimension greater than zero, there is a
procedure for obtaining an integer invariant as well: there is a tautological
complex line bundle $\Cal L$ on $M\times \Cal M(\xi)$, and one can use slant
product with $c_1(\Cal L)$ to obtain a class $\mu \in H^2(\Cal M(\xi); \Zee)$
whose top power can be integrated against the fundamental class of $\Cal
M(\xi)$ to give an integer. As long as
$b_2^+(M) > 1$, the invariants obtained are independent of the choice of the
Riemannian metric $g$, and in case $b_2^+(M) = 1$, there is a formula for the
dependence on the metric. 

Thus, assuming for simplicity that $b_2^+(M) > 1$ and that we have oriented 
$H^0(M; \Ar), H^1(M; \Ar)$, and $H^2_+(M; \Ar)$, we
have assigned an integer $SW(\xi)$ to every $\Spin ^c$ structure $\xi$ on $M$.
A compactness result shows that there are only finitely many $\xi$ such that
$SW(\xi) \neq 0$. If there exists a Riemannian metric $g$ on $M$ with   
positive scalar curvature, then in fact $SW(\xi) = 0$ for every $\Spin ^c$
structure 
$\xi$ on $M$. A $\Spin ^c$ structure $\xi$ such that $SW(\xi)
\neq 0$ will be called a {\sl basic\/} $\Spin ^c$ structure, and if  $L =\det
\xi$ where $\xi$ is a basic $\Spin ^c$ structure, then we will call $L$ a {\sl
basic class}. Witten has conjectured that, in case $M$ is of simple
type, then the basic classes in this sense are exactly the Kronheimer-Mrowka
basic classes $\kappa _i$, and has also conjecture the precise form of the
rational numbers
$a_i$ in the Kronheimer-Mrowka formula: up to a universal factor depending  
only on the homotopy type of $M$, they are just $\sum _{\det \xi = \kappa
_i}SW(\xi)$. Moreover, there is a corresponding notion of simple type in
Seiberg-Witten theory, that the expected dimension of the moduli space is 
always
$0$, and it is natural to expect that a $4$-manifold is of simple type in the
sense of Kronheimer-Mrowka if and only if it is of simple type in the sense of
Seiberg-Witten theory.

Let us make the following points about Seiberg-Witten theory:
\roster
\item Seiberg-Witten theory is essentially an abelian theory. Although the 
group
$\Spin ^c(4)$ is not abelian, the requirement that we only consider $\Spin
^c(4)$-bundles lifting the tangent bundle puts all of the extra degrees of
freedom in the choice of a complex line bundle on $M$. In particular the   
gauge group associated to Seiberg-Witten theory is abelian and the main
calculations involve the curvature of complex line bundles, as opposed to   
rank two complex vector bundles. Thus, in Yang-Mills theory, curvature
calculations do not give much useful information, whereas in Seiberg-Witten
theory they give the crucial compactness results as well as vanishing in the
case of positive scalar curvature.
\item The Seiberg-Witten equations are {\sl not\/} conformally invariant.
\item In connection with (1) and (2), no bubbling phenomena occur and the
moduli spaces are compact.
\item The fundamental group seems to play no major role in the definition or
the computation of the invariants.
\endroster
Thus the ideological underpinnings of Donaldson theory do not seem to address
the fundamental mystery of dimension $4$. Indeed, the Seiberg-Witten   
equations do not have anything like the natural interpretation of the 
Yang-Mills equations, and there does not as yet exist a mathematical
understanding of why the very complicated information of Donaldson theory 
should condense into the much more manageable information of Seiberg-Witten
theory. (A strategy for proving Witten's conjecture has been proposed by
Pidstrigach and Tyurin. This strategy gives a reason why there should be a  
link between the two theories, but it does not, at least to me, suggest why  
the Seiberg-Witten equations should be the right way to attack $4$-manifold
topology.)

\head 4. The case of a K\"ahler surface \endhead

In case $X$ is a K\"ahler surface, then $\Bbb S^+(\xi) = \Omega ^0(L_0)\oplus
\Omega ^{0,2}(L_0)$, we can write the Dirac operator in terms of the
$\dbar$-operator, and the Clifford multiplication used to identify the 
traceless endomorphisms from $\Bbb S^+(\xi)$ to itself with $\Omega ^2_+(M)
\otimes \Cee$ can be identified with wedge product or its adjoint on $\Omega
^0(L_0) \oplus \Omega ^{0,2}(L_0)$. A spinor field $\psi$ can be written
in terms of its components $(\alpha, \beta) \in \Omega ^0(L_0) \oplus
\Omega ^{0,2}(L_0)$. Working out the Seiberg-Witten
equations in this case gives the equations:
$$\align 
\dbar_A\alpha &+\dbar^*_A\beta =0; \\
F_A^{0,2} &= \dbar A^{0,1}
=\bar\alpha\beta ;\\
(F_A^+)^{1,1} &=\frac{i}{2}(|\alpha|^2-|\beta|^2)\omega .
\endalign$$
(Here the metric $g$ defines a Hermitian metric on $L_0$ and thus a conjugate
linear isomorphism $\Omega ^0(L_0) \to \Omega ^0(L_0^*)$, and $\bar \alpha$
denotes the image of $\alpha$ under this isomorphism.) We assume that
$(\psi,A)$ is an irreducible solution to the Seiberg-Witten equations. To
analyze the solutions, we argue as follows: applying $\dbar _A$ to the first
equation and plugging in the second, we get
$$ \dbar_A^2\alpha +\dbar _A\dbar^*_A\beta =0= F^{0,2}\alpha +\dbar
_A\dbar^*_A\beta = |\alpha|^2\beta + \dbar _A\dbar^*_A\beta .$$
Taking the inner product of this expression (namely $0$) with $\beta$ and
integrating over $X$ gives
$$\align
0 &= \int _X\left(\langle |\alpha|^2\beta, \beta \rangle + \langle\dbar
_A\dbar^*_A\beta , \beta \rangle\right) \\
&= \int _X |\alpha|^2 |\beta|^2 + \|\dbar^*_A\beta\|^2.
\endalign$$
Since both terms are positive, we must have $|\alpha|^2 |\beta|^2 = 0$
(pointwise) and $\|\dbar^*_A\beta\|^2 = 0$, so that $\dbar^*_A\beta = 0$.    
Now since $|\alpha|^2 |\beta|^2 = 0$, the product $\bar\alpha \beta = 0$ as 
well, so that $F_A^{0,2} = 0$. Thus as in the Yang-Mills case $A^{0,1}$ 
defines a holomorphic structure on $L$ and so on $L_0$. Since $\dbar^*_A  
\beta = 0$, the first equation implies that $\dbar_A\alpha = 0$, and so 
$\alpha$ is a holomorphic section of $L_0$ and $*\beta$ is a holomorphic 
section of
$K_X\otimes L_0^{-1}$. On the other hand, since $\bar\alpha \beta = 0$,  
$\beta$ must vanish on the open set $\{\alpha (x)\neq 0\}$. If this open set  
is nonempty, i\.e\. if $\alpha \neq 0$, then $\beta$ is identically zero, and
conversely if
$\beta \neq 0$ then $\alpha$ is identically zero. For simplicity let us   
assume that $\alpha \neq 0$. In this case, since $L_0$ has the nonzero
holomorphic section $L_0$, we can write $L_0 = \scrO_X(D_0)$ for an effective
curve $D_0$ on $X$. The first two equations for $(\alpha, \beta, A)$ just say
that $\beta = 0$, that $A$ defines a holomorphic structure on $L_0$, and that
$\alpha$ is a nonzero holomorphic section in the structure defined by $A$. In
other words, the $A$ and $\alpha$ give us effective curves $D_0$ on $X$ whose
associated cohomology class is equal to $c_1(L_0)$. Now the set of all curves 
on
$X$ with a fixed cohomology class is parametrized by a projective variety,   
the Hilbert scheme of $X$ (with an appropriate choice of Hilbert polynomial). 
The Hilbert scheme is an interesting and often highly nontrivial moduli space
associated to
$X$, and so as algebraic geometers we can get to work on this problem.

However, there is a third equation for $A$ and $\alpha$ which is part of the
Seiberg-Witten equations, namely the equation
$\dsize(F_A^+)^{1,1} =\frac{i}{2}(|\alpha|^2-|\beta|^2)\omega$. This equation 
for the
$(1,1)$ part of the curvature is analogous to the equation $\Lambda F_A^{1,1}
= 0$ in Yang-Mills theory, and should be thought of as a stability type
condition for the divisor $D_0$. As we shall see, however, the consequence of
stability in this case is strictly a numerical one. In case $\beta = 0$, the
equation becomes $\dsize(F_A^+)^{1,1} =\frac{i}{2}|\alpha|^2\omega$. In
particular since $c_1(L)$ is represented by the form $\dsize\frac{i}{2\pi
}F_A$, it follows that
$c_1(L)
\cdot \omega < 0$ in real cohomology. Thus we are led to the necessary
conditions for $(\alpha , 0, A)$ to be a solution to the Seiberg-Witten
equations:
\roster
\item $A$ defines a holomorphic structure on $L$;
\item In this structure, $\alpha$ defines a holomorphic section of $L_0 =
(K_X\otimes L)^{1/2}$;
\item $c_1(L) \cdot \omega < 0$.
\endroster
(In case $\alpha = 0, \beta\neq 0$, the appropriate change is: $\beta$   
defines a holomorphic section of $K_X\otimes L_0^{-1} = (K_X\otimes
L^{-1})^{1/2}$ and 
$c_1(L) \cdot \omega > 0$, and these conditions are Serre dual to the   
previous ones.) There is also the condition that the expected dimension of   
the moduli space is nonnegative, which for a complex surface is simply  
$L^2\geq K_X^2$, because $2\chi (X) + 3\sigma (X) = K_X^2$ by the Hodge index
theorem and Noether's formula. 

In fact, these necessary conditions are also sufficient. An easy calculation
with the K\"ahler identities shows that the existence of a solution involves
solving a PDE of the form
$$\Delta h - fe^{h/2} -g =0$$
for the unknown function $h$,
where $f$ and $g$ are given $C^\infty$ functions on $X$ such that $f\geq 0$
pointwise, where $f\neq 0$, and $\int _Xg < 0$, and $\Delta$ is the negative
definite Laplacian on $X$. This vortex equation has a unique $C^\infty$
solution, by a result of Kazdan-Warner \cite{20} (first exploited in gauge
theory by Bradlow \cite{4}). The analysis in the proof of this result is
nontrivial, but is much easier than that used in the proof of the
theorem of Donaldson and Uhlenbeck-Yau linking ASD connections and stable
bundles on a K\"ahler surface. In this way, the part of the Seiberg-Witten
moduli space corresponding to $\alpha \neq 0, \beta = 0$ (for the    
unperturbed equations) can be identified as a set with the Hilbert scheme of 
all effective curves $D_0$ on
$X$ such that, if we set $L =\scrO_X(2D_0) \otimes K_X^{-1}$, then $L$ 
satisfies the numerical condition
$c_1(L) \cdot \omega < 0$. There is a similar description in case $\alpha =
0,
\beta \neq 0$, which is Serre dual to the first one. Finally we should   
require that $L^2 \geq K_X^2$, since otherwise the moduli space will perturb
away to the empty set and the invariants we define will all be zero.

In this way we have identified the (unperturbed) Seiberg-Witten moduli space
with the union of finitely many components of the Hilbert scheme of curves    
on $X$. In fact, both spaces have additional structure: the Seiberg-Witten 
moduli space is locally modeled on a real analytic space, and the Hilbert 
scheme has a scheme structure and so is a complex analytic space (possibly
nonreduced). One can show that the two spaces are isomorphic as real analytic
spaces
\cite{17}. However, as we shall see in a minute, aside from a few very   
special cases this result is more of a virtual theorem than an actual method  
of computation.

We can now try to identify the Seiberg-Witten basic classes and then the  
moduli spaces for various algebraic surfaces. For an algebraic geometer, the 
end result is disappointing: the geometric interest of the basic classes and 
the moduli space is essentially inversely proportional to the interest of the
surface in question. Thus, if $X$ is a minimal surface of general type, the
basic classes are exactly $K_X^{\pm 1}$, corresponding to the trivial     
$\Spin ^c$ structures on $X$. Of course, it is exactly this fact which shows
that the classes $\pm[K_X]$ are preserved under diffeomorphisms. In fact, the
following is an easy consequence of the Hodge index theorem:

\proclaim{Theorem} Let $X$ be a minimal algebraic surface of general type,   
let $\omega$ be a K\"ahler form on $X$ and let $L$ be a holomorphic line  
bundle on $X$ such that:
\roster
\item $L^2\geq K_X^2$;
\item $\omega \cdot L < 0$;
\item There is an effective divisor $D_0$ such that $\scrO_X(2D_0)= K_X\otimes
L$.
\endroster
Then $L = K_X^{-1}$ and $D_0$ is the trivial divisor.
\endproclaim

\proclaim{Corollary} If $X$ is a minimal algebraic surface of general type,
then the only basic classes on $X$ are $\pm [K_X]$, with the trivial      
$\Spin ^c$ structures.
\endproclaim

\demo{Proof of the theorem} We write  $K_X$ and $L$ in the additive
notation of divisors, and will just consider the case corresponding to $\alpha
\neq 0$.  Since $\frac12(K_X+L)$ is effective, $\omega \cdot (K_X+L) \geq 0$, 
and likewise $K_X\cdot (K_X+L)
\geq 0$ since $K_X$ is nef. Thus $K_X\cdot L \geq -K_X^2$. On the other hand,
since $\omega \cdot L < 0$, there exists a $t\geq 1$ such that $\omega \cdot
(K_X+tL) =0$. By the Hodge index theorem, $(K_X+tL)^2 \leq 0$. Thus, since
$t\geq 1$ and $L^2\geq K_X^2$, 
$$\align
(K_X+tL)^2 &= K_X^2 + 2t(K_X\cdot L)+ t^2 L^2 \\
&\geq K_X^2 -2tK_X^2 + t^2K_X^2 = (1-2t+t^2)K_X^2 = (1-t)^2K_X^2 \geq 0.
\endalign$$
So $(K_X+tL)^2 = 0$, which can only happen if $t=1$ and 
$K_X+L$ is numerically trivial. Thus $D_0= \frac12(K_X+L)$ is numerically
trivial, and it is also effective, so that it is zero. It follows that the
corresponding line bundle
$L_0 = \scrO_X(D_0)$ is the trivial line bundle, and so $L= K_X^{-1}$ with the 
trivial $\Spin ^c$ structure.
\qed
\enddemo

We note that the value of the function $SW$ on $\pm [K_X]$ is $\pm 1$, which
follows formally since the trivial divisor is the unique smooth point of its
moduli space. In particular, $\pm [K_X]$ really are basic classes.

Similar arguments show the following:

\proclaim{Theorem} Let $X$ be a minimal algebraic surface with $\kappa
(X) \geq 0$ and $K_X^2 = 0$. Then the basic classes on $X$, as
elements of rational cohomology, are of the form $rK_X$ with $r\in \Bbb Q$   
and $|r| \leq 1$. Moreover the case $r=\pm 1$ does occur, and in this case   
the only corresponding $\Spin ^c$ structures are the trivial ones on        
$\pm [K_X]$.
\endproclaim

In the elliptic case, the invariant need not take on the value $\pm 1$. The
calculation of the value of the invariant is given in \cite{41} in case
$b_2^+\geq 3$ and in \cite{5} and \cite{17} in general. In case $X$ has at  
most two multiple fibers, the divisibility of $[K_X]$ and the largest $r\in 
\Bbb Q$ with $r\neq 1$ such that $r[K_X]$ is again a basic class will then
determine the possible multiplicities. In this way we can classify elliptic
surfaces up to diffeomorphism. 

For a nonminimal surface $X$, there is the following result, which is easy to
check directly from the holomorphic criteria for Seiberg-Witten classes:

\proclaim{Theorem} Let $X$ be a surface of general type with $b_2^+(X)      
\geq 3$, let $\rho \: X \to X_0$ be the blowdown to the minimal model, and   
let
$K_0$ be the preimage in $H^2(X; \Zee)$ of the canonical class of $X_0$ and
$E_1, \dots, E_r$ be the exceptional curves of $\rho$. Then the basic classes
for $X$ are exactly the classes $\pm K_0 \pm E_1\pm \cdots \pm E_r$, and the
invariant takes the value $\pm 1$ on each of these.
\endproclaim

Similar results hold for elliptic surfaces or in case $b_2^+= 1$, in which  
case the statement means that the classes above are the basic classes for a
particular chamber. It follows that, if $b_2^+(X) \geq 3$, then $X$ is
always of simple type in the Seiberg-Witten sense. Using the above theorem,  
one can show that the classes
$K_0$ and the $E_i$ are $C^\infty$ invariants up to sign and permutation of  
the
$E_i$. In fact this is clear up to $2$-torsion for the case $b_2^+\geq 3$,   
and follows by an analysis of the chamber structure in case $b_2^+=1$  
\cite{5},
\cite{16}. In case $b_2^+(X) \geq 3$ and $X$ is minimal, we see that
the Seiberg-Witten basic classes are equal to $\pm [K_X]$ if $X$ is of   
general type, and are of the form $r[K_X]$ where $r\in \Bbb Q$ in general. 
Thus, if Witten's conjecture on the form of the Donaldson polynomials is true,
then the Donaldson polynomials of a minimal algebraic surface $X$ with 
$b_2^+(X)
\geq 3$ always lie in $\Cee[q_X, k]$, where $k=[K_X]$, and actually involve 
$k$, as long as $k$ is nonzero in rational cohomology.

There remains the case where $X$ is rational or ruled. In this case   
$b_2^+(X) = 1$, and the invariants depend on a chamber structure. It is easy  
to see from the holomorphic criterion for Seiberg-Witten invariants or from  
the existence of metrics with positive scalar curvature on $X$ that there is
always a chamber on $X$ for which all of the invariants vanish, i\.e\. for 
which there are no basic classes. On the other hand, there are chambers where 
the invariants are nonzero. Some of these values correspond to moduli spaces
where the expected dimension of the Seiberg-Witten moduli space is nonzero,
i\.e\. $X$ does not necessarily have simple type in case $b_2^+(X) = 1$. For 
the case of rational surfaces the nonzero values are always
$\pm 1$. However, for irrational ruled surfaces the values are more
complicated, and the associated moduli spaces involve interesting algebraic
geometry. For example, in case $X = C\times \Pee ^1$, the moduli spaces are
connected with correspondences on $X$ and tie in with the theory of special
divisors on $C$. In case $X$ is a general ruled surface over $C$, the
structure of the Seiberg-Witten moduli space involves the study of the   
Hilbert scheme on $X$ and is related to questions concerning general rank two
stable vector bundles over $C$. For a discussion of these results, see
\cite{17}. However, while the algebraic geometry involved is nontrivial,
Seiberg-Witten theory in this case does not seem to have any consequences for
the $C^\infty$ topology of $X$, which can be analyzed by elementary methods.

\head 5. Concluding remarks \endhead

The ability of the Seiberg-Witten invariants to solve seemingly intractable
questions on the smooth topology of algebraic surfaces is a stunning
achievement. Of course, like all such achievements, there is also a certain
amount of disappointment: the Seiberg-Witten invariants, and therefore
presumably also Donaldson invariants, are only able to tell us about the
pullback of the canonical class of the minimal model and the exceptional  
curves (at least for $\kappa (X) \geq 0$), and these classes are the obvious
ones. Is this the end of the story as far as gauge theory invariants are
concerned? There are surfaces of general type, for example some of the  
Horikawa surfaces, which are homeomorphic but not deformation equivalent, and
which cannot be distinguished by Seiberg-Witten or Donaldson invariants. 
Perhaps some new multi-monopole invariants will be able to show that such
surfaces are not diffeomorphic, perhaps completely new invariants are needed, 
or perhaps the surfaces are indeed diffeomorphic but not deformation 
equivalent. Still other questions, such as the problem of understanding the
mapping class groups of algebraic surfaces or more general
$4$-manifolds, or in other words understanding the difference between
pseudo-isotopy and isotopy in dimension $4$, remain and do not seem   
accessible either. Beyond these questions, the unruly world of surfaces of
general type, for which no reasonable classification is known to exist, is   
now known to account for only a small part of the even more complicated world 
of smooth
$4$-manifolds. 

On the positive side, even as the story on gauge theory and algebraic   
surfaces seems to be coming to a conclusion, or at least a natural pause,    
the study of symplectic $4$-manifolds has opened up dramatically thanks mainly 
to deep new work of Taubes \cite{39}. Whereas for algebraic surfaces the
classification theory leads to a simple description of the Seiberg-Witten
invariants, in the case of symplectic $4$-manifolds it is the invariants
themselves which have led to a deep series of results paralleling the
classification theory of algebraic surfaces. This story is still in progress.

\refstyle{A}


\Refs

\ref \no  1\by R. Barlow\paper A simply connected surface of general type
with $p_g=0$\jour Inventiones Math.\vol 79\yr 1985\pages 293--301\endref

\ref \no  2\by S. Bauer\paper Some nonreduced moduli of bundles and
Donaldson invariants for Dolgachev surfaces\jour J. reine angew. Math.\vol
424\yr 1992\pages 149--180\endref

\ref \no  3\bysame \paper  Diffeomorphism classification of elliptic
surfaces with $p_g=1$\jour J. reine angew. Math.\vol
451\yr 1994\pages 89--148 \endref


\ref \no 4 \by S. Bradlow \paper Vortices for holomorphic line bundles
over closed K\"ahler manifolds \jour Comm. Math. Physics \vol 135 \yr
1990 \pages 1--17\endref

\ref \no 5 \by R. Brussee \paper Some $C^\infty$ properties of K\"ahler
surfaces \paperinfo Algebraic geometry e-prints 9503004 \endref

\ref \no  6\by I. Dolgachev \paper Algebraic surfaces with $q = p_g = 0$
\inbook in Algebraic Surfaces \bookinfo C.I.M.E. Cortona 1977 \publ Liguori
\publaddr Napoli  \yr 1981 \pages 97--215 \endref

\ref \no 7\by S.K. Donaldson \paper Anti-self-dual Yang-Mills connections  
over  complex algebraic surfaces and stable vector bundles\jour Proc. Lond.
Math.  Soc. \vol 50\pages 1--26 \yr 1985\endref

\ref \no 8\bysame \paper Irrationality and the $h$-cobordism
conjecture\jour J. Differential Geom. \vol 26\pages 141--168 \yr 1987\endref


\ref \no 9\bysame\paper Polynomial invariants for smooth 
four-manifolds \jour Topology \vol 29 \pages  257--315 \yr 1990\endref

\ref \no 10\bysame\paper The Seiberg-Witten equations and
$4$-manifold topology \jour Bull. Amer. Math. Soc.
(N.S.) 
\vol 33 \pages  45--71 \yr 1996\endref


\ref \no 11\by S.K. Donaldson and P. Kronheimer \book The Geometry of 
Four-Manifolds
\publ Clarendon Press \publaddr Oxford \yr 1990 \endref

\ref \no 12\by R. Friedman
\paper Vector bundles and $SO(3)$-invariants for elliptic surfaces \jour   
Jour. Amer. Math. Soc. \vol 8 \yr 1995 \pages 29--139
\endref

\ref \no 13\by R. Friedman and  J. W.  Morgan\paper On the
diffeomorphism types of certain algebraic surfaces I \jour J. Differential 
Geom.
\vol   27
\pages  297--369 \yr 1988 \moreref \paper  II \jour J. Differential
Geom. \vol 27 \yr 1988 \pages 371--398 \endref

\ref \no 14\bysame \paper Algebraic surfaces and
$4$-manifolds:  some conjectures and speculations \jour Bull. Amer. Math.   
Soc. (N.S.) 
\vol   18 \pages  1--19 \yr 1988\endref

\ref \no 15\bysame\book Smooth Four-Manifolds and
Complex Surfaces, {\rm Ergebnisse der Mathematik und ihrer Grenz\-gebiete 3. 
Folge} {\bf 27} \publ Springer \publaddr Berlin Heidelberg  New York \yr
1994\endref

\ref \no 16\bysame \paper Algebraic surfaces and Seiberg-Witten invariants
\toappear \endref

\ref \no 17\bysame \paper Obstruction bundles, semiregularity, and
Seiberg-Witten invariants \jour Comm. Analysis and Geometry \toappear
\endref

\ref \no 18\by R. Friedman and Z.B. Qin \paper On complex surfaces 
diffeomorphic to rational surfaces \jour Inventiones Math. \vol 120 \yr 1995
\pages 81--117
\endref

\ref \no  19\by J. Harer, A. Kas and R. Kirby \book Handlebody
Decompositions of Complex Surfaces \bookinfo Memoirs of the American
Mathematical Society {\bf 350} \publ American Mathematical Society
\publaddr Providence \yr 1986 \endref


\ref \no 20 \by J. Kazdan and F. W. Warner \paper Curvature functions for
compact $2$-manifolds \jour Annals of Math. \vol 99 \yr 1974 \pages
14--47 \endref

\ref \no  21\by K. Kodaira \paper On homotopy $K3$ surfaces \inbook  in Essays 
on Topology and Related Topics \bookinfo M\'emoires d\'edi\'es \`a Georges de
Rham \publ Springer \yr 1970 \pages 58--69 \endref

\ref \no  22\by P. Kronheimer \paper Embedded surfaces in $4$-manifolds  
\inbook Proceedings of the International Congress of Mathematicians Kyoto   
1990
\publ The Mathematical Society of Japan/Springer Verlag \publaddr Tokyo   
Berlin Heidelberg New York \yr 1991 \pages 529--539 \endref

\ref \no  23\by P. Kronheimer and T. Mrowka \paper Gauge theory for embedded
surfaces I\jour Topology \vol 32 \yr 1993 \pages 773--826 \moreref \paper II 
\jour Topology \vol 34 \yr 1995 \pages 37--97
\endref


\ref \no  24\bysame \paper Recurrence relations and
asymptotics for four-manifold invariants \jour Bull. Amer. Math. Soc.
(NS) \vol 30 \yr 1994 \pages 215--221 \endref

\ref \no 25  \bysame \paper Embedded surfaces and the structure of
Donaldson's polynomial invariants \jour J. Differential Geometry \vol 41
\yr 1995 \pages 573--734\endref 

\ref \no 26  \bysame \paper The genus of embedded surfaces in the projective
plane \jour Math. Research Letters \vol 1 \yr 1994 \pages 797--808 \endref

\ref \no  27\by B. Moishezon \book Complex Surfaces and Connected Sums of
Complex Projective Planes \bookinfo Lecture Notes in Mathematics {\bf 603}
\publ Springer Verlag \publaddr Berlin Heidelberg New York \yr 1977 \endref

\ref \no  28\by J. W. Morgan \book The Seiberg-Witten Equations and
Applications to the Topology of Smooth Four-Manifolds \bookinfo  Mathematical
Notes
\vol 44 \publ Princeton University Press
\publaddr Princeton \yr 1996
\endref

\ref \no  29\by J. W. Morgan and T. Mrowka \paper A note on Donaldson's
polynomial invariants 
\jour International Math. Research Notices  \yr 1992 \pages 223--230
\endref

\ref \no  30\bysame \paper On the diffeomorphism
classification of regular elliptic surfaces \jour International Math.
Research Notices \vol 6 \yr 1993 \pages 183--184 \endref

\ref \no  31\by J. W. Morgan and K. O'Grady \book Differential Topology of
Complex Surfaces Elliptic Surfaces with $p_g=1$: Smooth Classification
\bookinfo Lecture Notes in Mathematics \vol 1545 \publ Springer-Verlag
\publaddr Berlin Heidelberg New York \yr 1993
\endref

\ref \no 32 \by J. W. Morgan, Z. Szab\'o, and C. H. Taubes \paper A product
formula for Seiberg-Witten invariants and the generalized Thom conjecture
\toappear
\endref

\ref \no 33 \by C. Okonek and A. Teleman \paper The coupled Seiberg-Witten
equations, vortices, and moduli spaces of stable pairs \paperinfo Algebraic
geometry e-prints 9505012
\endref

\ref \no 34 \bysame \paper Seiberg-Witten invariants and the Van de Ven
conjecture
\paperinfo Algebraic geometry e-prints 9505013
\endref

\ref \no 35 \by V.Y. Pidstrigach \paper Patching formulas for spin
polynomials, and a proof of the Van de Ven conjecture \jour Russian Academy  
of  Science Izvestiya Mathematics, Translations of the AMS \vol 45 \yr 1995
\pages 529--544 \endref

\ref \no 36\by V.Y. Pidstrigach and A.N.  Tyurin \paper Invariants
of the smooth  structure of an algebraic surface arising from the Dirac 
operator
\jour  Russian Academy of Science Izvestiya Mathematics, Translations of the
AMS 
\vol   40 \pages  267--351 \yr 1993\endref

\ref \no 37 \by N. Seiberg and E. Witten \paper Electric-magnetic duality,
monopole condensation, and confinement  in $N=2$ supersymmetric Yang-Mills
theory \jour Nuclear Physics B \vol 426 \yr 1994 \pages 19--52 \endref

\ref \no 38 \by F. Severi \paper La g\'eom\'etrie alg\'ebrique italienne. Sa
rigueur, ses m\'ethodes, ses probl\`emes \inbook Colloque de g\'eom\'etrie
alg\'ebrique, Li\`ege 1949 \publ Georges Thone \publaddr Li\`ege \moreref  
\publ Masson et Cie \publaddr Paris \yr 1950 \pages 9--55 \endref

\ref \no 39 \by C. H. Taubes \paper Seiberg-Witten invariants and Gromov
invariants \toappear \endref

\ref \no 40 \by K. Uhlenbeck and S.-T. Yau \paper On the existence of
Hermitian Yang-Mills connections on stable vector bundles \jour Comm. Pure
Appl. Math. \vol 39 \yr 1986 \pages 257--293 \moreref \paper A note on our
previous paper: On the existence of
Hermitian Yang-Mills connections on stable vector bundles \jour Comm. Pure
Appl. Math. \vol 42 \yr 1989 \pages 703--707 \endref

\ref \no 41\by E. Witten \paper Monopoles and four-manifolds \jour Math.
Research Letters \vol 1 \yr 1994 \pages 769--796
\endref

\ref \no 42 \by O. Zariski \book Algebraic Surfaces \rom(Second Supplemented
Edition\rom) \bookinfo Ergebnisse der Mathematik und ihrer Grenzgebiete     
\vol 61
\publ Springer Verlag \publaddr New York Heidelberg Berlin \yr 1971 \endref

\endRefs
\enddocument